\documentclass[12pt]{amsart}

\usepackage{fullpage, graphicx, url, xcolor}

\newtheorem{thm}{Theorem}[section]
\newtheorem{prop}[thm]{Proposition}

\newtheorem{lemma}[thm]{Lemma}

\newtheorem{definition}[thm]{Definition}
\theoremstyle{definition}

\title{A study on prefixes of $c_2$ invariants}

\author{Karen Yeats}

\thanks{Thanks to Oliver Schnetz for goading me into doing the prefix calculation, and for his continued excitement about the $c_2$ invariant.  Thanks to Iain Crump and Freddy Cachazo for discussions.  Thanks to NSERC and the Humboldt foundation for support.  Thanks to the organizers of the Algebraic Combinatorics, Resurgence, Moulds and Applications (CARMA) conference for their excellent conference.}

\begin{document}
\maketitle

\begin{abstract}
  This paper begins by reviewing recent progress that has been made by taking a combinatorial perspective on the $c_2$ invariant, an arithmetic graph invariant with connections to Feynman integrals.
  Then it proceeds to report on some recent calculations of $c_2$ invariants for two families of circulant graphs at small primes.  These calculations support the idea that all possible finite sequences appear as initial segments of $c_2$ invariants, in contrast to their apparent sparsity on small graphs.
\end{abstract}

\section{Introduction}

The $c_2$ invariant is an arithmetic graph invariant introduced by Schnetz in \cite{SFq} in order to better understand certain Feynman integrals.  The $c_2$ invariant sees aspects of the same underlying geometry that the Feynman period sees \cite{BrS, BrS3} and consequently the $c_2$ invariant can predict things about what classes of numbers can show up in a given Feynman period.  See Subsection~\ref{subsec CARMA} for some further comments in this direction.  In the following graphs will be assumed to be connected unless otherwise mentioned.
For a graph $G$, the $c_2$ invariant of $G$, $c_2(G)$ is a sequence of numbers indexed by primes
\[
c_2(G) = (c_2^{(2)}(G), c_2^{(3)}(G), c_2^{(5)}(G), c_2^{(7)}(G), c_2^{(11)}(G), \ldots) 
\]
with each $c_2^{(p)}(G) \in \mathbb{Z}/p\mathbb{Z}$.  The definition of $c_2^{(p)}(G)$ is given in the next section.  As we will see, $c_2^{(q)}$ can be defined for prime powers, not just primes, but we will stick to primes herein.  
Previous work of Brown and Schnetz \cite{BrS3} calculated $c_2$ invariants for graphs up to 10 loops at small primes (up to the first 100 primes in some cases).  These calculations uncovered many interesting patterns, most notably coefficient sequences of $q$-expansions of modular forms.

In \cite{Ycirc}, the author described a new, more graphical technique for calculating $c_2$ invariants.  With Wesley Chorney, this technique was expanded in \cite{CYgrid}.  For this technique the prime $p$ is fixed but there is a finite algorithm to calculate the $c_2$ invariant at $p$ for all members of a recursively constructed class of graphs.  The loop orders of the graphs in these classes are unbounded, so these techniques let us calculate $c_2$ invariants at all loop orders albeit with fixed $p$ and only for certain families.  

The $c_2$ invariant is also believed to have certain symmetries corresponding to symmetries of the Feynman period.  One of these is called completion symmetry and was conjectured by Brown and Schnetz in 2010 in \cite{BrS}.  This conjecture has turned out to be quite difficult.
This combinatorial perspective on the $c_2$ invariant is used in \cite{Yscompl} to prove one special case of the conjecture.  An overview of the results of \cite{Ycirc}, \cite{CYgrid}, and \cite{Yscompl} is given below in Section~\ref{sec past results} along with an outlook for this approach and connections to topics of particular interest to the CARMA conference.

The explicit calculations in \cite{Ycirc} and \cite{CYgrid} were done by hand and so involve only $p=2$ and relatively simple classes of graphs.  The complexity of the calculation as a function of $p$ is quite bad, and it also grows depending on the graph class.  Nonetheless, by computer some new progress is possible, though only a little, the results of which are reported on in Section~\ref{sec comp}.  These new computations are particularly interesting because they let us probe the behaviour of the $c_2$ invariant at all loop orders, giving a rather different impression than previous exhaustive computations at fixed loop orders.

\section{Set up}\label{sec setup}

Let $G$ be a 4-regular graph.  The graph resulting from removing any one vertex of $G$ (and its adjacent edges) is called a \emph{decompletion} of $G$, and $G$ is called the \emph{completion} of any of its decompletions.  In general there may be many non-isomorphic decompletions of a graph.

For any graph $H$ (but of primary interest is the case when $H$ is a decompletion of a 4-regular graph) associate an indeterminate $a_e$ for each edge $e$ and define the (dual) \emph{Kirchhoff polynomial} or \emph{first Symanzik polynomial} of $G$ to be
\[
\Psi_H = \sum_{T}\prod_{e\not\in T}a_e
\]
where the sum runs over all spanning trees $T$ of $H$.  For example the Kirchhoff polynomial of a 3-cycle with edge variables $a, b, c$ is $a+b+c$.

Now we can define the $c_2$ invariant.  Given a polynomial $f$ with integer coefficients, write $[f]_q$ for the number of $\mathbb{F}_q$-rational points on the affine variety defined by $f=0$ (with $f$ first reduced to $\mathbb{F}_q$).  Our polynomials $f$ will always come from a graph in one way or another, and so the affine space in which they are to be taken will always be of dimension the number of edges of the graph.

\begin{definition}
  Suppose $H$ has at least 3 vertices, then 
  \[
  c_2^{(q)}(H) = \frac{[\Psi_H]_q}{q^2} \mod q.
  \]
\end{definition}
See \cite{BrS} for a proof that this is well-defined.  Note that in \cite{BrS} they have the condition that the dimension of the cycle space of $H$ is at most two less than the number of edges of $H$, however using Euler's formula this condition is equivalent to $H$ having at least 3 vertices.  In what follows we will restrict to $c_2$ invariants at primes $p$, not more general prime powers, since this is computationally accessible and corresponds to what has been calculated elsewhere \cite{BrS3}.

We can also view $\Psi_H$ as a determinant in the following way. Choose an arbitrary order for the edges and the vertices of $H$ and choose an arbitrary orientation for the edges of $H$.  Let $E$ be the signed incidence matrix of $H$ with one row removed and let $\Lambda$ be the diagonal matrix of the edge variables of $H$.  Let
\[
M = \begin{bmatrix} \Lambda & E^t \\ -E & 0 \end{bmatrix}.
\]
Then by the matrix-tree theorem 
\[
\det(M) = \Psi_H.
\]
(See \cite{Brbig} Proposition 21 or \cite{VY} for details.)  The matrix $M$
behaves much like the Laplacian matrix of a graph with variables included 
and with one matching row and column removed, but the pieces which make it up are expanded out by blocks, so call $M$ the \emph{expanded Laplacian} of $H$.

As well as $\det(M)$, minors of $M$ are useful.  For $I$ and $J$ sets of edge indices (or sets of edges; with the edge order fixed, we need not distinguish between an edge and its index) let $M(I,J)$ be the expanded Laplacian with rows indexed by elements of $I$ removed and columns indexed by elements of $J$ removed. In \cite{Brbig} Brown gave the following definition.
\begin{definition}
Let $I$, $J$, and $K$ be sets of edge indices with $|I|=|J|$.  Define
\[
\Psi^{I,J}_{G,K} = \det(M(I,J))|_{\substack{a_e=0 \\ \text{for $e\in K$}}}
\]
When $K=\emptyset$ we will simply leave it out.
\end{definition}
Brown called these polynomials \emph{Dodgson polynomials}.  They satisfy many relations, see \cite{Brbig}. Different choices in the construction of $M$
may change the overall sign of a Dodgson polynomial, but since we will be concerned with counting zeros of these polynomials the overall sign is of no interest.

The combinatorial perspective on the $c_2$ invariant comes from now taking a different view point.  Instead of thinking about the $c_2$ invariant in terms of polynomials and point counts, we want to think about it in terms of set partitions of subsets of vertices and spanning forest polynomials.  We need a few lemmas to get to this reinterpretation.

\begin{lemma}[Lemma 24 of \cite{BrS} along with inclusion-exclusion]\label{lem D3}
  Suppose $2 + |E(H)| \leq 2|V(H)|$.  Let $i,j,k$ be distinct edge indices of $H$ and let $p$ be a prime.  Then
  \[
  c_2^{(p)}(H) = -[\Psi^{ik,jk}_{H}\Psi^{i,j}_{H,k}]_p \mod p.
  \]
\end{lemma}
Note that if $H$ is a decompletion of $4$-regular graph with at least 2 vertices then it satisfies the hypotheses of the previous lemma.  This lemma is useful because we no longer have to divide by $p^2$ but rather directly count points modulo $p$.  Combined with the next lemma, we no longer need to count points at all.

\begin{lemma}\label{lem cor of CW}
  Let $F$ be a polynomial of degree $N$ in $N$ variables, $x_1, \ldots, x_N$, with integer coefficients. The coefficient of $x_1^{p-1}\cdots x_N^{p-1}$ in $F^{p-1}$ is $[F]_p$ modulo $p$.
\end{lemma}
This lemma is a corollary of one of the standard proofs of the Chevalley-Warning theorem, see section 2 of \cite{Ax}.

Together these two lemmas tell us that to calculate the $c_2$ invariant we only need to understand the coefficient of 
\[
\prod_{\substack{1 \leq \ell \leq |E(H)| \\ \ell \neq i,j,k}}a_\ell^{p-1} \quad  \text{in}  \quad \left(\Psi^{ik,jk}_{H}\Psi^{i,j}_{H,k}\right)^{p-1}.
\]  We can make this yet more combinatorial by reinterpreting these Dodgson polynomials in terms of sums over spanning forests.  

\begin{definition}
  Let $P$ be a set partition of a subset of the vertices of $H$.  Define
  \[
  \Phi^P_H = \sum_{F}\prod_{e \not\in F}a_e
  \]
  where the sum runs over spanning forests $F$ of $H$ with a bijection between the trees of $F$ and
the parts of $P$
where each vertex in a part lies in its corresponding tree.  Trees consisting of isolated vertices are allowed.
\end{definition}
Call these polynomials \emph{spanning forest polynomials}.
Dodgson polynomials can always be rewritten in terms of spanning forest polynomials.  This is a manifestation of the all-minors matrix tree theorem \cite{Chai}; the following form is convenient for the present purposes.
\begin{prop}[Proposition  12  from  \cite{BrY}]
Let $I$, $J$, and $K$ be sets of edge indices with $|I|=|J|$.  Then
\[
\Psi^{I,J}_{H,K} = \sum\pm \Phi^P_{H\backslash(I\cup J\cup K)}
\]
where the sum runs over all set partitions $P$ of the end points of edges of $(I \cup J\cup K) \backslash (I\cap J)$ such that all the forests corresponding to $P$ become spanning trees in both $G\backslash I/(J\cup K)$ and $G\backslash J/(I\cup K)$.
\end{prop}
The signs in the sum can be determined, see Proposition 16 of \cite{BrY}.  All that we will need is that if the set partitions are of the form $\{a,b\},\{c,d\}$ and $\{a,c\}, \{b,d\}$ then they appear with opposite sign, see Corollary 17 of \cite{BrY}.

The two lemmas told us to calculate the coefficient of $\prod_{\ell \neq i,j,k}a_\ell^{p-1}$ in $\left(\Psi^{ik,jk}_{H}\Psi^{i,j}_{H,k}\right)^{p-1}$ modulo $p$.  Now, we can interpret the two Dodgson polynomials as signed sums of spanning forest polynomials, and so we are interested in the coefficient of $\prod_{\ell \neq i,j,k}a_\ell^{p-1}$ in each of certain products of $2p-2$ spanning forest polynomials.  Summing those coefficients and taking the result modulo $p$ then calculates the $c_2$ invariant.  Notice that each spanning forest polynomial is, by construction, linear in each edge variable.  So taking this coefficient amounts to determining which of the variables to assign to each polynomial in the product and taking the resulting monomial from each polynomial.  If a particular variable is assigned to a particular spanning forest polynomial then we are restricting ourselves to the spanning forests in that polynomial which do not use that edge.  If a particular variable is not assigned to a particular spanning forest polynomial then we are restricting ourselves to the spanning forests in that polynomial which do use that edge, or equivalently, to spanning forests in the graph with that edge removed with one more tree than before made from breaking up the tree which originally used that edge.

Given a subset $S$ of the edges of a graph $H$ and a spanning forest polynomial $\Phi^P_H$, an \emph{assignment} of the edges of $S$ to $H$ is the polynomial resulting from a choice for each edge of $S$ to either assign it or not assign it to $\Phi^P_H$.  That is, given a choice of a subset $S'\subseteq S$, the resulting polynomial is the coefficient of $\prod_{e\in S'}a_e$ in $\Phi^P_H|_{a_e=0, e\in S-S'}$.  This new polynomial is itself a sum of spanning forest polynomials, as the following lemma describes.

\begin{lemma}[Lemma 5.2 from \cite{CYgrid}]\label{lem edge red}
Given a spanning forest polynomial on a graph $H$ and a set $S\subseteq E(H)$, any assignment of the edges of $S$ yields a sum of spanning forest polynomials on the graph $H-S$. Furthermore, the vertices involved in the set partitions defining the new spanning forest polynomials involve only vertices already in partition for the original polynomial along with vertices incident to $S$.
\end{lemma}
Note that some set partitions may give impossibilities, in this case the spanning forest polynomial is an empty sum, and so is $0$.  Also, we can discard isolated vertices from $H-S$ as by connectivity they must each be in their own part of every set partition and so they contribute no information.

\section{Past applications of this method}\label{sec past results}

\subsection{Completion}

One of the important reasons to study the $c_2$ invariant is to better understand the Feynman period.  For a graph $H$ the Feynman period, in affine form, is defined to be
\[
P_H = \int_{a_i \geq 0} \frac{\prod da_e}{\Psi_H^2}\bigg|_{a_1=1}.
\]
This is a residue of the Feynman integral in parametric form which is independent of kinematical parameters.  If $H$ comes from a four regular graph $K$ with one vertex removed and $K$ is internally 6-edge-connected, that is any way of removing fewer than six edges of $K$ either leaves the graph connected or disconnects only an isolated vertex, then the integral converges.  The Feynman period is known to have four important symmetries, most of which were long known in physics, and which can be found in a form as we will use them in \cite{Sphi4}.   The symmetry we will focus on is the completion symmetry, namely if $H_1$ and $H_2$ are both decompletions of the same 4-regular graph $G$ then $P_{H_1} = P_{H_2}$.  This is proved for the period by moving to momentum space and inverting the variables, see \cite{Sphi4}.

The Feynman period is controlled by the geometry of the denominator of the integral, that is by the geometry of the variety $\Psi_H=0$.  The $c_2$ invariant is accessing this geometry from a different direction, by counting rational points in that same variety over various finite fields.  Thus we should expect that they are saying something about each other.  However, the $c_2$ invariant is only seeing part of the geometric structure, not the whole thing.  As it turns out the $c_2$ invariant has been very useful in predicting properties of the period, most notably if $c_2^{(p)}=0$ for all $p$ then we expect the period to have less than maximal transcendental weight for the size of the graph.  We also expect that if two graphs have the same period then they should have the same $c_2$ invariant.  The converse is certainly not true, with $c_2^{(p)}=0$ for all $p$ being a good example.  One consequence of this is that the $c_2$ invariant should have all the symmetries that the period does.

In particular Brown and Schnetz, \cite{BrS}, conjectured that if $G$ is a connected 4-regular graph and $v$ and $w$ are vertices of $G$ then $c_2^{(p)}(G-v) = c_2^{(p)}(G-w)$ for all primes $p$.

The main result of \cite{Yscompl} is a very special case of this conjecture and the first major progress towards the conjecture.
\begin{thm}[Theorem 1.2 of \cite{Yscompl}]
  Let $G$ be a connected 4-regular graph with an odd number of vertices.
  Let $v$ and $w$ be vertices of $G$.  Then $c_2^{(2)}(G-v) = c_2^{(2)}(G-w)$. 
\end{thm}

The approach of \cite{Yscompl} is combinatorial following the set up described in Section~\ref{sec setup}.  For $p=2$ this is particularly simple since we are looking for the coefficient of $\prod_{\ell \neq i,j,k}a_\ell$ in $\Psi^{ik,jk}_{H}\Psi^{i,j}_{H,k}$ modulo $2$; that is, we need to determine the parity of the number of edge assignments compatible with $\Psi^{ik,jk}_{H}\Psi^{i,j}_{H,k}$, and each edge assignment assigns exactly one copy of each edge, dividing them between the two polynomials, so we are counting certain edge bipartitions.

Further, it suffices to prove the result for $v$ and $w$ adjacent, leaving the remaining vertex between $v$ and $w$ 3-valent in each decompletion, and we can take $i,j,k$ to be the three incident edges to this vertex.  This means that for both $G-v$ and $G-w$ we are down to considering spanning forest polynomials on the graph $G-\{v,w\}$.  Further, then, in either decompletion $\Psi^{ik, jk} = \Psi_{G-\{v,w\}}$, and so we only need to count the parity of the number of edge bipartitions such that one part is a spanning tree of $G-\{v,w\}$ and the other part is a spanning forest compatible with $\Psi^{i, j}_k$; which spanning forests these are is the only thing that changes between $G-v$ and $G-w$.

What we need is for the number of these edge bipartitions to have the same parity between $G-v$ and $G-w$.  This takes some work; the bipartitions fall into two classes depending on how the vertex partitions giving the spanning forest polynomials divide the vertices.  For some of the cases, we can define a fixed point free involution showing that the set of these edge bipartitions is even.  For other cases, no such construction was evident, and so instead a more complicated construction was used involving an auxiliary graph related to the spanning tree graph.  What was needed to finish the proof was that this auxiliary graph has an even number of vertices, as the vertices correspond to the edge bipartitions in the remaining cases.  When $G$ has an odd number of vertices, all the vertices of the auxiliary graph have odd degree and so the auxiliary graph as a whole has an even number of vertices.  When $G$ has an even number of vertices then we do not have this parity restriction on the vertex degrees in the auxiliary graph, and so the proof does not extend directly.

\subsection{Circulants and toroidal grids}\label{subsec circ}

Another use of the approach described in Section~\ref{sec setup} is to calculate $c_2$ invariants for fixed $p$ but for whole families of graphs.  This contrasts with \cite{BrS3} where Brown and Schnetz fix the graph and calculate the $c_2$ invariant for many primes.  They do this systematically for small graphs, collecting many interesting results.

Taking the fixed $p$ and graph family approach, in \cite{Ycirc} and \cite{CYgrid} we have both specific and general results.  For the specific results we need to define some classes of graphs.
\begin{definition} The circulant graph $C_n(i_1, i_2, \ldots, i_k)$ is the graph on $n$ vertices with an edge between vertices $i$ and $j$ if and only if $i-j = i_\ell \mod n$ or $j-i = i_\ell \mod n$ for some $1\leq \ell \leq k$.
\end{definition}
We will be interested in certain $4$-regular circulant graphs.  Every $4$-regular circulant graph can be written as $C_n(i,j)$ with $i\neq n-j$ and $i,j\neq n/2$.  Note that each decompletion of a given circulant graph is isomorphic, so we will use the notation $\widetilde{G}$ for the decompletion of $G$, which is well-defined in the case that $G$ is a circulant or other vertex-transitive graph.

Circulant graphs are an interesting class of graphs for questions related to $c_2$ or to Feynman periods because they include both the simplest non-trivial family, namely the zigzags \cite{Szigzag}, which are the $C_n(1,2)$, but they also include some of the most difficult and mysterious graphs; for example consider the last few entries at each loop order in \cite{Sphi4}.  So circulant graphs cut across difficulties while being graph theoretically very nice because of their symmetries.  For the present purposes, what is even better is that they are very well suited to the algorithm we will discuss below.

\begin{definition}
The nonskew toroidal grid indexed by $(k, 0)$ and $(0, m)$ is the Cartesian product
of a cycle on $k$ vertices and a cycle on $m$ vertices.
\end{definition}

These are nonskew toroidal grids as we can also define more general toroidal grids where we begin with the Cartesian product of a cycle on $k$ vertices and a path with $m$ edges, making a finite cylindrical grid, and then identify the top and bottom $k$-cycles, potentially with an offset.  In the case where there is a nonzero offset $\ell$ then the result is a skew toroidal grid and if $\gcd(m,\ell)=1$ then it is isomorphic to the circulant graph $C_{km}(\ell, m)$

The author (\cite{Ycirc}) and the author with Wesley Chorney (\cite{CYgrid}) proved the following explicit results
\begin{thm}
\mbox{}
  \begin{itemize}
    \item $c_2^{(2)}(\widetilde{C_n}(1,3)) = n \mod 2$ for $n\geq 7$ (\cite{Ycirc} Proposition 4.1).
    \item $c_2^{(2)}(\widetilde{C_{2k+2}}(1,k)) = 0 \mod 2$ for $k\geq 3$ (\cite{Ycirc} Proposition 7.1).
    \item Let $G$ be a nonskew toroidal grid indexed by $(N, 0)$, $(0, m)$ with at least $3N$ vertices.  Then $c_2^{(2)}(\widetilde{G}) = 0$ (\cite{CYgrid} Proposition 3.2)
  \end{itemize}
\end{thm}
In \cite{CYgrid} we also prove that two other families have $c_2^{(2)}=0$.

The technique behind these results is an algorithm which applies to any \emph{recursively constructible family} of graphs with the property that Lemma~\ref{lem D3} eventually applies to members of the family.  The notion of recursively constructible family is due to Noy and Rib\'o, see \cite{NRrec} section 2.  For its use in this $c_2$ context see \cite{CYgrid} section 5.  We do not need the precise definition here, but the idea is as follows.  

Roughly, recursively constructible families consist of graphs made with some initial piece and then a chain of repeated structures and then a cap which may link the last piece of the chain back to the initial piece.  Let $\{H_n\}$ be the graphs of this family with $n$ the length of the chain.  Let $H'_n$ be $H_n$ with the edges of the cap deleted.   Fix a prime $p$.  Using Lemma~\ref{lem D3} with three edges in the cap and using Lemma~\ref{lem edge red} to assign all other edges in the cap, we can calculate $c_2^{(p)}(H_n)$ by taking the coefficient of $\prod_{e\in H'_n} a_e^{p-1}$ in some sum of products of $2p-2$ spanning forest polynomials of $H'_n$ where the partitions only use vertices in the final piece of the chain and in the initial piece.  Note that there are only a finite number of spanning forest polynomials of this form and so only a finite number of products of $2p-2$ spanning forest polynomials of this form.  For each product of spanning forests appearing, we can further use Lemma~\ref{lem edge red} to assign all edges of one piece of the chain.  Removing one piece of chain gives us $H'_{n-1}$ so what we obtain from Lemma~\ref{lem edge red} is a sum of products of $2p-2$ spanning forest polynomials of $H'_{n-1}$.  What this tells us is that we can calculate the coefficient of $\prod_{e\in H'_n} a_e^{p-1}$ in any of these products as a sum of coefficient of $\prod_{e\in H'_{n-1}} a_e^{p-1}$ in products for $H'_{n-1}$.  So what we have is a system of first order linear recurrences.  The system itself can be computed in a finite amount of time (at worst all possible products of $2p-2$ spanning forest polynomials need to be reduced in this way).  For the initial conditions we need to compute the required coefficient for each product of polynomials for $H'_{1}$.  Then it remains to solve the system modulo $p$.  This algorithm is due to the author with Wesley Chorney and gives the following theorem
\begin{thm}[\cite{CYgrid}, Theorem 5.3]
Let ${G_n}$ be a recursively constructible family of graphs with $2|V(G_n)| =|E(G_n)|+2$ for $n$ sufficiently large. The $c_2$ invariant for any fixed prime $p$ can be calculated
using these methods in a finite amount of time for all graphs of the family.
\end{thm}
The special case of this algorithm for the circulant graphs used in Section~\ref{sec comp} was already in \cite{Ycirc}.

As a piece of mathematics the existence of this algorithm is very nice as it says that it is theoretically possible to rigorously compute $c_2$ invariants for entire families of graphs, albeit only for fixed $p$.  In practice it is not so nice.  The specific results listed above were obtained by applying this method on the specific examples using hand computations.  Note that we stuck to the particularly nice case of $p=2$ for the explicit results.  The complexity of the method as a function of $p$ is very poor.  The complexity of the method also grows quickly in the number of vertices which $H_{n-1}'$ shares with the next piece of chain. 

\subsection{Comments and outlook from past results}\label{subsec CARMA}

The $c_2$ invariant began as a tool linking the Feynman period, which is in some sense a physical object, to the arithmetic of the Kirchhoff variety.  These come from different aspects of the geometry of the Kirchhoff variety.  This is why some arithmetic questions on the kinds of numbers appearing in the Feynman integrals and their transcendental properties have a fairly tight link to the $c_2$ invariant and its properties.  With the methods surveyed above another major perspective and toolset, that of combinatorics, is brought to the study of the $c_2$ invariant and progress can be made.

The completion result, partial though it is, is a testament to the power of algebraic and enumerative combinatorics.  By doing some counting using classic enumerative tools like fixed-point-free involutions progress was made where more algebro-geometric tools were stuck.  In the end both areas are important as the reduction to the counting problem is fundamentally arithmetic in nature.  This interplay between areas is part of what makes these problems enjoyable and is central to the aims of the CARMA conference and the broader CARMA project.

We can also ask: what does it mean that all nonskew toroidal grids have $c_2^{(2)}=0$?  When the $c_2$-invariant is $0$ for all $p$ then we expect a drop in transcendental weight; it is too much to hope that this is what is happening here, though to clarify this one of the most interesting explicit calculations to do now would be $p=3$ for some nonskew toroidal grids.  More likely the symmetries of the toroidal grids only force that $c_2^{(2)}=0$.  It is not clear what this means geometrically, nor, in the other direction, what other graphs behave in this way, though some other examples are known.

Towards the future, the obstacles to extending the partial completion result do not appear insurmountable and are the subject of ongoing research.  Getting beyond $p=2$ will involve seeing how the larger number of possibilities collect into sets of size $p$.  We can also try to collect more data from the family approach.  At this point in the study of the $c_2$ invariant more data nearly invariably shows new patterns and raises new questions.  The remainder of this document is a report on a computerization of the circulant $c_2$ calculation in the $\widetilde{C_n}(1,3)$ and $\widetilde{C_n}(2,3)$ cases and a discussion of some of the perhaps surprising things which can be seen in this data.

\section{Computerized circulant $c_2$ computations}\label{sec comp}

\begin{figure}
\includegraphics{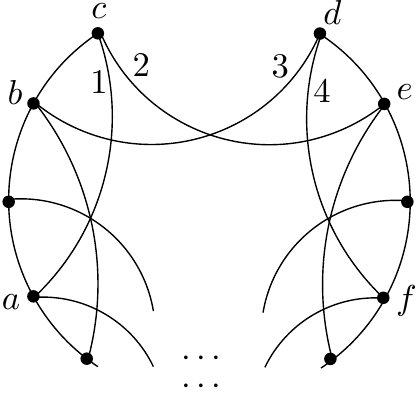} \qquad \includegraphics{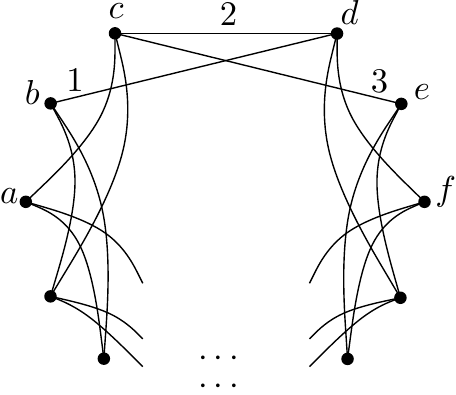}
\caption{$\widetilde{C_n}(1,3)$ (left) and $\widetilde{C_n}(2,3)$ (right).}\label{fig our graphs}
\end{figure}

The goal now is to implement the algorithm described in Subsection~\ref{subsec circ} in as practical a manner as possible on certain families of circulant graphs.



From now on we will only consider the families $\widetilde{C_n}(1,3)$ and $\widetilde{C_n}(2,3)$.   These graphs are illustrated in Figure~\ref{fig our graphs}.  The first step of the algorithm is to process enough edges so that what remains has a chain structure. 

For $\widetilde{C_n}(1,3)$ we process the edges labelled 1,2,3, and 4 in Figure~\ref{fig our graphs} and for $\widetilde{C_n}(2,3)$ we process the edges labelled 1,2, and 3.  This can be done by hand with the lemmas of the previous section, and in fact these calculations were done in \cite{Ycirc}.  As it turns out, for $\widetilde{C_n}(1,3)$ dealing with the other edges incident to $c$ and $d$ comes for free, giving
\begin{align*}
  c_2^{(p)}(\widetilde{C_n}(1,3)) & = [\prod_{e\in H_n}a_e^{p-1}]\left(\Phi_{H_n}^{\{a,f\}, \{b\}, \{e\}}\Psi_{H_n}\right)^{p-1} \mod p, \\
  c_2^{(p)}(\widetilde{C_n}(2,3)) & = [\prod_{e\in K_n}a_e^{p-1}]\left(\Phi_{K_n}^{\{b,e\},\{c\},\{d\}}\left(\Phi_{H_n}^{\{c,d\}, \{b,e\}}-\Phi_{H_n}^{\{c,b\}, \{d,e\}} \right)\right)^{p-1} \mod p
\end{align*}
where here square brackets are the combinatorialist's notation for \emph{coefficient of}\footnote{Specifically, if $F$ is a polynomial and $m$ is a monomial then $[m]F$ is the coefficient of $m$ in $F$.}, where $H_n$ is $\widetilde{C_n}(1,3)$ with vertices $c$ and $d$ and their incident edges deleted, and where $K_n$ is $\widetilde{C_n}(2,3)$ with edges $1$, $2$, and $3$ deleted.  The spanning forest reductions can be found in Section~4 of \cite{Ycirc} for $C_n(1,3)$ and Section~5 of \cite{Ycirc} for $C_n(2,3)$. 

The next step of the algorithm is to build the system of linear recurrences by processing the edges of one piece of the chain for each product of spanning forest polynomials which can appear.  In these cases, one piece of chain is the single end vertex and its incident edges, say $e$ for $H_n$ and $d$ for $K_n$.  Note that for the products of spanning forest polynomials calculated in the previous paragraph, at most the first three and last three vertices of the graph are used in the set partitions.  In view of Lemma~\ref{lem edge red}, this remains true (now in $H_{n-1}$ or $K_{n-1}$) after assigning the edges of one piece of the chain.  Thus the products of spanning forest polynomials which could appear are any products of $2p-2$ spanning forest polynomials where each polynomial in the product comes from a partition of the first 3 and last 3 vertices and where the total degree of the product is correct.  Throughout this algorithm, we never explicitly deal with polynomials.  A spanning forest polynomial is always represented by its partition of these 6 vertices and all the manipulations work directly with the partitions.

There are 203 set partitions of 6 elements (\cite{OEIS} A000110), but fortunately, nowhere near all of possible products of them are actually required.  The implementation makes a list of products which need to be processed, starting with the ones which give $c_2$ itself and adding to the list as needed.  The number of products, call this $N$, which were necessary for the program in each case computed is shown in Table~\ref{tab N}.  For the actual edge assignments, we are only reducing two edges at a time, so for each polynomial in the product there are only four possibilities, both edges are in, both edges are out, the one edge alone is in, or the other edge alone is in.  How each of these affects the vertex partition is coded for each case and then the edge assignment calculation simply comes down to looping over the ways of assigning edges between the polynomials.  The outcome of this step of the algorithm is the $N\times N$ coefficient matrix of the system of recurrences.

\begin{table}
\begin{tabular}{|lr|}
\hline
$C(1,3)$ & \\
\hline
$p$ & $N$ \\
\hline
2 & 29 \\
3 & 546 \\
5 & 82703 \\
7 & 5698505\\
\hline
\hline
$C(2,3)$ & \\
\hline
$p$ & $N$ \\
\hline
2 & 248 \\
3 & 30729 \\
\hline
\end{tabular}
\caption{Number of products of spanning forest polynomials ($N$) necessary.}\label{tab N}
\end{table}

The next step of the algorithm is to calculate the initial conditions.  We need $N$ initial conditions.  Each of these is calculated on the smallest graph of the family.  The smallest $H_n$ has 5 vertices and 6 edges; the smallest $K_n$ has 5 vertices and 5 edges.  This is the case where the last of the first three vertices is the same as the first of the last three vertices.  For both the minimal $H_n$ and minimal $K_n$ cases, how each possible spanning forest partitions the vertices is precomputed, then for a given product of spanning forest polynomials each partition is compared to the precomputed list to get the count.  The outcome of this step is a vector of length $N$.

The final step is then to solve this system of linear recurrences with these initial conditions.  This is not done algebraically because the matrices get very large.  Rather, the system is simply iterated.
The sequence of $c_2$ invariants at $p$ is then the sequence of first entries of these iterated vectors in the $H_n$ case and a weighted sum of the first $p+1$ entries in the $K_n$ case.  
Eventually these sequences seem to begin to repeat.  It turns out that the vector at the point where the $c_2$ sequence first repeats is not yet equal to the initial condition vector.  Rather, it takes multiple iterations of the repeating block of the $c_2$ sequence before the vector matches the initial condition vector.  Call the period of repetition of the $c_2$ sequence the \emph{$c_2$ period} and call the period of repetition of the vector the \emph{vector period}\footnote{Note the different use of the term period from earlier sections: for the remainder of the paper, period will mean period of repetition of a sequence}.  Verifying that the vector agrees with the initial condition vector at the vector period proves that the system repeats with this period and then checking the $c_2$ sequence breaks into blocks according to the $c_2$ period within one vector period proves that the $c_2$ sequence repeats with the observed $c_2$ period as well.

Unfortunately, the vector periods are quite large making them computationally problematic.  The easiest and most naive way to compute the vector period is simply to iterate the system and compare the resulting vector with the initial condition vector.  This process is not guaranteed to terminate as there could be transient behaviour in the early iterations.  The simplest example of such transients would be if there was a row of all $0$s but a nonzero initial condition in that location, but longer transients are also possible.  All we are guaranteed theoretically is that, by finiteness of the field, at some point the result of an iteration agrees with some past iteration.  In principle this is also true of the $c_2$ itself, but in practice the $c_2$ sequence displays periodic behaviour beginning at the very first value.  A less naive way to compute the vector period would compare it with past values after each $c_2$ period.  The downside of this approach is that it is slower and uses more memory.  

Attempts were made to compute the vector periods for $p=5$ for $C(1,3)$ and for $p=3$ for $C(2,3)$.  For $p=5$ and $C(1,3)$ the less naive computation method was used but ran out of memory after three weeks; the vector period in this case exceeds 153844320.  For $p=3$ and $C(2,3)$ the more naive computation run for a month suggests that the vector period exceeds 4614354360, though it remains possible that initial transient behaviour simply means that the vector period is not obtainable by the naive method.  An attempt to compute the vector period for $p=7$ for $C(1,3)$ was not made because iterating the system until the point where the $c_2$ sequence appeared to repeat took over a month and the vector period would be expected to be many times this.  

The periods are given for each computed case in Table~\ref{tab periods}.  Note that the $c_2$ periods are not proved, only empirically observed, for $p=5$ and $p=7$ in the $C(1,3)$ case, nor for $p=3$ for the $C(2,3)$ case.  However, the empirical evidence is quite strong.  For $p=5$ and $C(1,3)$, the block of the $c_2$ sequence repeated exactly 41356 times before the computation was killed.  For $p=3$ and $C(2,3)$, the block of the $c_2$ sequence repeated exactly 1059310 times before this document was submitted.

\begin{table}
\begin{tabular}{|lrr|}
\hline
$C(1,3)$ & & \\
\hline
$p$ & $c_2$ period & vector period\\
\hline
2 & 2 & 4\\
3 & 36 & 59040\\
5 & 3720 &  \\
7 & 134064 &  \\
\hline
\hline
$C(2,3)$ & & \\
\hline
$p$ & $c_2$ period & vector period\\
\hline
2 & 7 & 56 \\
3 & 4356 &  \\
\hline
\end{tabular}
\caption{Periods of repetition for the system of recurrences.  The $c_2$ periods are proven when the vector period is listed and are otherwise empirical.}\label{tab periods}
\end{table}

Observe that the $c_2$ periods are all much smaller than one would naively expect given the sizes of the matrices and also considerably smaller than the vector periods.  This means that there is a substantial amount of structure which this method does not capture.  The ratio between the $c_2$ periods and the vector periods gives the first hint of where this additional structure may reside.  Looking at the vectors after each $c_2$ period, the first many entries agree while some later entries do not.  The system is built so that for each entry of the vector the corresponding product of spanning forest polynomials does appear in the construction, but the behaviour of the vector and the $c_2$ sequence indicates that various values of the later entries are equivalent for the $c_2$ calculation.  Playing around with the coefficient matrix in the $C(1,3)$, $p=2$ case, suggests that some block decomposition might be possible in order to explain at least some of the discrepancy between the $c_2$ period and the vector period.  Unfortunately, the structure is not clear for the $p=3$ coefficient matrix.  Understanding this redundancy should be the next step for both the theoretical and practical take on this algorithm.

Without such additional reductions, the computations presented here exhaust what we can do for $\widetilde{C_n}(1,3)$ and $\widetilde{C_n}(2,3)$.  The $p=7$ computation for $\widetilde{C_n}(1,3)$ took 100GB of RAM and took several months to run on a University of Waterloo server.  Even then the system was only iterated until the $c_2$ sequence appeared to repeat.  Specifically, the system was iterated until the first 1351 entries repeated and prior to this there was no reoccurrence of an initial segment of length greater than 6.  Consequently, $p=11$ will be outside the range of practical computation.  The $p=5$ case for $\widetilde{C_n}(2,3)$ was attempted but was killed as it exceeded 400GB of RAM; $N$ had already surpassed 10 million and rough heuristics based on how $N$ grew during the other computations suggests that the final $N$ for $p=5$ for $\widetilde{C_n}(2,3)$ is likely to be in excess of 100 million.

Finally, then, the $c_2$ invariants themselves are presented.  All the sequences have been verified for small values of $n$ by Oliver Schnetz using different techniques.

\begin{align*}
  c_2^{(2)}(\widetilde{C_n}(1,3)) & = (1,0)^* \\
  c_2^{(3)}(\widetilde{C_n}(1,3)) & = (0, 0, 0, 0, 0, 0, 1, 2, 2, 1, 2, 2, 2, 2, 1, 1, 1, 2, 0, 1, 0, 2, 0, 1, 1, 1, 2, 2, 2, 1, 2, 0, 1, 0, 1, 0)^* \\
  c_2^{(2)}(\widetilde{C_n}(2,3)) & = (1, 1, 1, 0, 1, 0, 0)^*
\end{align*}
The remaining computed $c_2$ sequences, $c_2^{(5)}(\widetilde{C_n}(1,3))$, $c_2^{(7)}(\widetilde{C_n}(1,3))$, and $c_2^{(3)}(\widetilde{C_n}(2,3))$ are included in the auxiliary files in the arXiv submission as is the code.

The sequences as presented for $\widetilde{C_n}(1,3)$ begin at $n=9$ which corresponds to 7 loop decompleted graphs.  The sequences for $\widetilde{C_n}(2,3)$ begin at $n=7$ which corresponds to 5 loop decompleted graphs.  The $*$ indicates to repeat the sequence indefinitely.  For example $c_2^{(2)}(\widetilde{C_n}(1,3)) = (1,0)^*$ means that $c_2^{(2)}(\widetilde{C_9}(1,3)) = 1$, $c_2^{(2)}(\widetilde{C_{10}}(1,3)) = 0$, $c_2^{(2)}(\widetilde{C_{11}}(1,3)) = 1$ and so on.  Note that $c_2^{(2)}(\widetilde{C_n}(1,3)$ was computed by hand in \cite{Ycirc} and that the three sequences displayed above are proved while the remaining three are only empirically observed.

\section{Discussion}

This data is interesting and important because it lets us probe $c_2$ invariants at all loop orders, albeit only on these two families of graphs and only for a very few initial primes.

One particularly interesting question is which finite sequences occur as initial sequences of $c_2$ invariants $(c_2^{(2)}(H), c_2^{(3)}(H), c_2^{(5)}(H), c_2^{(7)}(H))$ in this data.  This is easy to tally, we just take the least common multiple of the periods to get the period for the initial segments so far and then count how many of each occur.  To begin with consider the prefix $(c_2^{(2)}, c_2^{(3)})$.  There are 6 possible prefixes and the distributions for each family are shown in Table~\ref{tab 2 prefix}.  Note that for $\widetilde{C_n}(1,3)$ it is uniform, while for $\widetilde{C_n}(2,3)$ the difference between the largest and smallest counts is less than $6\%$ of the total number.

\begin{table}
\begin{tabular}{|l|r|r|}
\hline
prefix & $C(1,3)$& $C(2,3)$ \\
\hline
 & count  & count  \\
\hline
(0,0) &  6 & 4236  \\
(0,1) &  6 & 4389  \\
(0,2) &  6 & 4443  \\
(1,0) &  6 & 5648  \\
(1,1) &  6 & 5852  \\
(1,2) &  6 & 5924  \\
\hline
\end{tabular}
\caption{Frequencies for prefixes $(c_2^{(2)}, c_2^{(3)})$.}\label{tab 2 prefix}
\end{table}

\begin{figure}
  \scalebox{0.8}{
\setlength{\unitlength}{1pt}
\begin{picture}(0,0)
\includegraphics{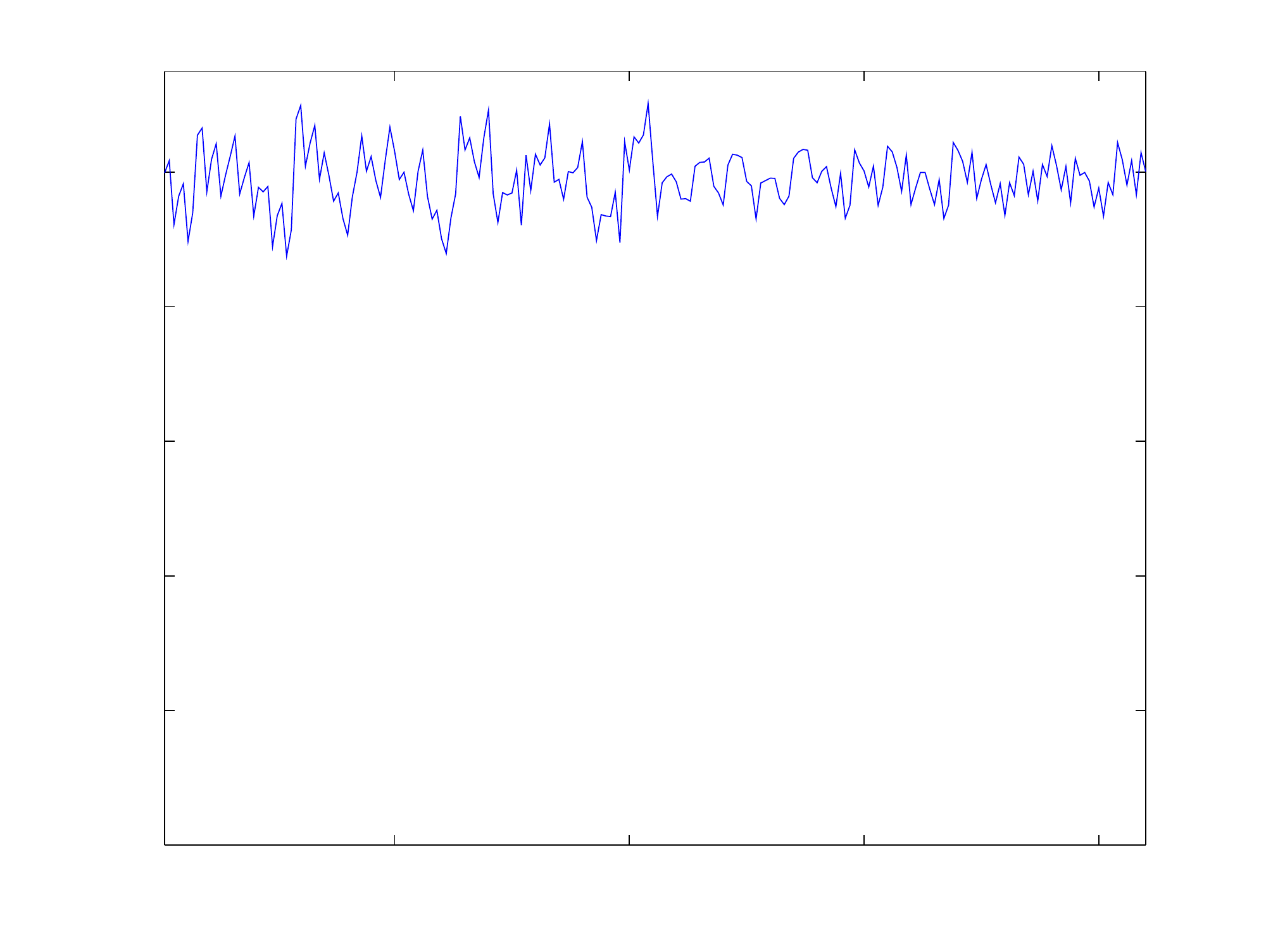}
\end{picture}%
\begin{picture}(576,433)(0,0)
\fontsize{10}{0}
\selectfont\put(179.538,43.5189){\makebox(0,0)[t]{\textcolor[rgb]{0,0,0}{{50}}}}
\fontsize{10}{0}
\selectfont\put(286.333,43.5189){\makebox(0,0)[t]{\textcolor[rgb]{0,0,0}{{100}}}}
\fontsize{10}{0}
\selectfont\put(393.127,43.5189){\makebox(0,0)[t]{\textcolor[rgb]{0,0,0}{{150}}}}
\fontsize{10}{0}
\selectfont\put(499.921,43.5189){\makebox(0,0)[t]{\textcolor[rgb]{0,0,0}{{200}}}}
\fontsize{10}{0}
\selectfont\put(69.8755,48.52){\makebox(0,0)[r]{\textcolor[rgb]{0,0,0}{{0}}}}
\fontsize{10}{0}
\selectfont\put(69.8755,109.751){\makebox(0,0)[r]{\textcolor[rgb]{0,0,0}{{20000}}}}
\fontsize{10}{0}
\selectfont\put(69.8755,170.983){\makebox(0,0)[r]{\textcolor[rgb]{0,0,0}{{40000}}}}
\fontsize{10}{0}
\selectfont\put(69.8755,232.214){\makebox(0,0)[r]{\textcolor[rgb]{0,0,0}{{60000}}}}
\fontsize{10}{0}
\selectfont\put(69.8755,293.445){\makebox(0,0)[r]{\textcolor[rgb]{0,0,0}{{80000}}}}
\fontsize{10}{0}
\selectfont\put(69.8755,354.677){\makebox(0,0)[r]{\textcolor[rgb]{0,0,0}{{100000}}}}
\fontsize{10}{0}
\selectfont\put(298.08,32.5189){\makebox(0,0)[t]{\textcolor[rgb]{0,0,0}{{prefix (index in lexicographic order)}}}}
\fontsize{10}{0}
\selectfont\put(29.8755,224.56){\rotatebox{90}{\makebox(0,0)[b]{\textcolor[rgb]{0,0,0}{{count}}}}}
\end{picture}}
\caption{Number of occurrences of length 4 prefixes}\label{fig plot}
\end{figure}
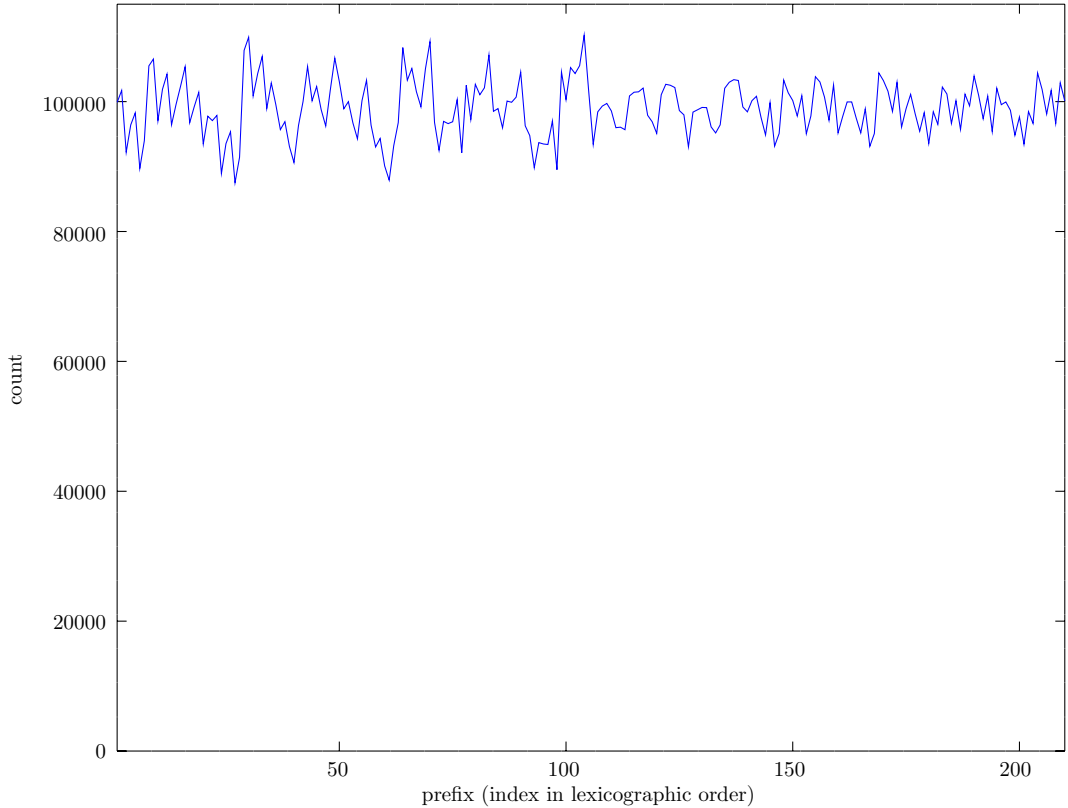

For $\widetilde{C_n}(1,3)$ we can also consider the prefixes of length 3 and 4.  For the prefixes of length 3 the period for the prefix is 11160.  The numbers of occurrences of the prefixes are between 350 and 393; the difference between these is less than $0.4\%$ of the total.  The mean is 372.  Performing the same calculations on the prefixes of length 4.  The period for the prefix is 20779920.  The counts all lie between 87514 and 110213; the difference between these is slightly over $0.1\%$ of the total.  The mean of the counts is 98952.  The number of occurrences for the prefixes of length 4 are plotted in Figure~\ref{fig plot}.

{}From this we see that the distribution of the frequencies of the different prefixes is quite flat and does not seem to be getting any less flat as we take larger prefixes.  In particular every prefix occurs and there is no indication that this will change as we move to longer prefixes.  This is an interesting and perhaps unexpected observation because looking at small graphs leaves the impression that only rather few finite sequences occur as prefixes of $c_2$ invariants.  This data, which can probe all loop orders, suggests quite the opposite: perhaps \emph{all} finite sequences can occur as prefixes of $c_2$ invariants.  

Geometrically, this says that the possible geometries for Kirchhoff varieties should not be expected to be sparse among all possible geometries, rather, at least as far as finite prefixes can see, it looks like everything can happen, though it may take very large loop order to get there.

\bibliographystyle{plain}
\bibliography{main}

\end{document}